\documentstyle[12pt,amscd,amssymb,verbatim]{amsart}
\numberwithin{equation}{section}
\renewcommand{\rm}{\normalshape}%

\theoremstyle{plain}
\newtheorem{theorem}{Theorem}
\newtheorem{prop}[theorem]{Proposition}
\newtheorem{proposition}[theorem]{Proposition}

\newtheorem{corollary}[theorem]{Corollary}

\theoremstyle{definition}

\theoremstyle{remark}
\newtheorem{remark}[theorem]{Remark}
\newtheorem{example}[theorem]{Example}

 \def\today{\ifcase\month\or
  January\or February\or March\or April\or May\or June\or
  July\or August\or September\or October\or November\or December\fi
  \space\number\day, \number\year}

\def\cA{{\cal A}}
\def\cB{{\cal B}}
\def\cC{{\cal C}}
\def\cD{{\cal D}}

\def\cQ{{\cal Q}}

\def\Br{\textrm{Br}}

\def\tr{\operatorname {tr}}

\begin{document}
\title[Genus of division algebras]
{On genus of division algebras}
\author[Sergey V. Tikhonov ]
{Sergey V. Tikhonov}

\address{
Belarusian State University, Nezavisimosti Ave., 4,
220030, Minsk, Belarus} \email{tikhonovsv@@bsu.by }

\begin{abstract}
The genus ${\bf gen}(\cD)$ of a finite-dimensional central division algebra $\cD$ over a field $F$ is defined as the collection of classes $[\cD']\in \Br(F)$,
where $\cD'$ is a central division $F$-algebra having the same maximal subfields as $\cD$.
We show that the fact that quaternion division algebras $\cD$ and $\cD'$ have the same maximal subfields does not imply that the matrix algebras $M_l(\cD)$ and $M_l(\cD')$ have
the same maximal subfields for $l>1$. Moreover, for any odd $n>1$, we construct a field $L$ such that there are two quaternion division $L$-algebras
$\cD$ and $\cD'$ and a central division $L$-algebra $\cC$ of degree and exponent $n$ such that ${\bf gen} (\cD) = {\bf gen} (\cD')$ but
${\bf gen} (\cD \otimes \cC) \ne {\bf gen} (\cD' \otimes \cC)$.

\end{abstract}

\maketitle

\def\dd{{\partial}}

\def\toeq{{@>\sim>>}}
\def\into{{\hookrightarrow}}

\def\emptyset{{\varnothing}}

\def\alp{{\alpha}}  \def\bet{{\beta}} \def\gam{{\gamma}}
 \def\del{{\delta}}
\def\eps{{\varepsilon}}
\def\kap{{\kappa}}                   \def\Chi{\text{X}}
\def\lam{{\lambda}}
 \def\sig{{\sigma}}  \def\vphi{{\varphi}} \def\om{{\omega}}
\def\Gam{{\Gamma}}  \def\Del{{\Delta}}  \def\Sig{{\Sigma}}
\def\ups{{\upsilon}}


\def\A{{\mathbb A}}
\def\F{{\mathbb F}}
\def\Q{{{\mathbb{Q}}}}
\def\CC{{\mathbb{C}}}
\def\PP{{\mathbb P}}
\def\R{{\mathbb R}}
\def\Z{{\mathbb Z}}
\def\X{{\mathbb X}}
\def\N{{\mathbb N}}
\def\C{{\mathbb C}}

\def\Gm{{{\Bbb G}_m}}
\def\Gmk{{{\Bbb G}_{m,k}}}
\def\GmL{{\Bbb G_{{\rm m},L}}}
\def\Ga{{{\Bbb G}_a}}

\def\Fb{{\overline F}}
\def\Hb{{\overline H}}
\def\Kb{{\overline K}}
\def\Lb{{\overline L}}
\def\Yb{{\overline Y}}
\def\Xb{{\overline X}}
\def\Tb{{\overline T}}
\def\Bb{{\overline B}}
\def\Gb{{\overline G}}
\def\Vb{{\overline V}}

\def\kb{{\bar k}}
\def\xb{{\bar x}}

\def\Th{{\hat T}}
\def\Bh{{\hat B}}
\def\Gh{{\hat G}}

\def\Xt{{\tilde X}}
\def\Gt{{\tilde G}}

\def\gg{{\mathfrak g}}
\def\gm{{\mathfrak m}}
\def\gp{{\mathfrak p}}
\def\gq{{\mathfrak q}}

\def\min{^{-1}}

\def\textrm#1{\text{\rm #1}}

\def\char{\textrm{char}}
\def\cor{\textrm{cor}}
\def\tr{{\rm{tr}}}
\def\res{{\rm{res}}}
\def\ram{{\rm{ram}}}
\def\deg{\textrm{deg}}
\def\Nrd{\textrm{Nrd}}
\def\N{\textrm{N}}
\def\exp{\textrm{exp}}
\def\Gal{\textrm{Gal}}
\def\Spec{\textrm{Spec}}
\def\Proj{\textrm{Proj}}
\def\Perm{\textrm{Perm}}
\def\coker{\textrm{coker\,}}
\def\Hom{\textrm{Hom}}
\def\im{\textrm{im\,}}
\def\ind{\textrm{ind}}
\def\int{\textrm{int}}
\def\inv{\textrm{inv}}

\def\tors{_{\textrm{tors}}}      \def\tor{^{\textrm{tor}}}
\def\red{^{\textrm{red}}}         \def\nt{^{\textrm{ssu}}}
\def\sc{^{\textrm{sc}}}
\def\sss{^{\textrm{ss}}}          \def\uu{^{\textrm{u}}}
\def\ad{^{\textrm{ad}}}           \def\mm{^{\textrm{m}}}
\def\tm{^\times}                  \def\mult{^{\textrm{mult}}}
\def\tt{^{\textrm{t}}}
\def\uss{^{\textrm{ssu}}}         \def\ssu{^{\textrm{ssu}}}
\def\cf{^{\textrm{cf}}}
\def\ab{_{\textrm{ab}}}

\def\et{_{\textrm{\'et}}}
\def\nr{_{\textrm{nr}}}

\def\op{^{\textrm{op}}}

\def\til{\;\widetilde{}\;}

\def\emptyset{{\varnothing}}



The genus ${\bf gen}(\cD)$ of a finite-dimensional central division algebra $\cD$ over a field $F$ is defined as the collection of classes $[\cD']\in \Br(F)$,
where $\cD'$ is a central division $F$-algebra having the same maximal subfields as $\cD$.
This means that $\cD$ and $\cD'$ have the same degree $n$, and a field extension $K/F$ of degree $n$ admits an $F$-embedding $K \hookrightarrow \cD$ if and only if it admits an $F$-embedding
$K \hookrightarrow \cD'$. Different variations of the notion of the genus are mentioned in \cite{ChRaRa15}.

The following questions were formulated in \cite[footnote 1 and Remark 2.2]{ChRaRa18}:

\medskip

{\it Does the fact that division algebras $\cD$ and $\cD'$ have the same maximal subfields imply that the matrix algebras $M_l(\cD)$ and $M_l(\cD')$ have
the same maximal subfields / \'{e}tale subalgebras for any (or even some) $l>1$?
}

\medskip

{\it Let $n_1$ and $n_2$ be relatively prime positive integers. Let also $\cD_i$ and $\cD'_i$ be central division algebras of degree $n_i$ over a field $F$ for $i=1,2$.
Is it true that if ${\bf gen}(\cD_i)={\bf gen}(\cD'_i)$ for $i=1,2$, then ${\bf gen}(\cD_1\otimes\cD_2)={\bf gen}(\cD'_1\otimes\cD'_2)$?
}

\medskip

Negative answers to these questions are given in Theorem \ref{th: main} and Corollary \ref{cor:matrix algebras} below.

We use the following notation.
For a field $F$, $F^*$ denotes the multiplicative group of $F$. ${F^*}^2$ denotes the subgroup of squares in $F^*$.
For a field extension $K/F$ and a central simple $F$-algebra $\cA$, $\cA_K$ denotes the tensor product
$\cA\otimes_F K$ and $\res_{K/F} : \Br(F)\longrightarrow \Br(K)$ denotes the restriction homomorphism.

Let $F$ be a field, $\char(F) \ne 2$. Assume that there are two non-isomorphic quaternion division $F$-algebras $\cA$ and $\cB$.
We will be interested in finite separable field extensions $K/F$ that satisfy the following three conditions:


(A) There is $d \in K^* \backslash {K^*}^2$ such that
there is no $a\in F^*$ such that
$da \in {K^*}^2$;

(B) $K$ does not split $\cB$;

(C) $K(\sqrt{d})$ splits $\cB$ but does not split $\cA$.

%


\begin{example} \label{ex:field}
Let $M$ be a field containing a primitive $n$th root of 1, $\char(M) \nmid 2n$. Let also
$F := M(x,y,z,w)$ be a purely transcendental extension of $M$ of transcendence degree 4 and $\cA := (x,y)$, $\cB := (w,y)$ quaternion $F$-algebras.
Then $K := F(\sqrt[n]{z})$ is a cyclic extension of $F$ of degree $n$ and $K$ does not split $\cB$. Finally, let $d := w+y (\sqrt[n]{z}+1)^2$.
Since $K(\sqrt{d})$ is a purely transcendental extension of $M(x,y)$, then $K(\sqrt{d})$ does not split $\cA$.
On the other side, $K(\sqrt{d})$ splits $\cB$ since $d$ is represented over $K$ by the quadratic form $<w,y,-yw>$. Finally, note that there is no $a\in F^*$ such that
$da \in {K^*}^2$. 
Thus the extension $K/F$ satisfies conditions (A) - (C).
Note also that the symbol $F$-algebra $(x,z)_n$ of degree and exponent $n$ is split by $K$.
Moreover, if  $\char(M)=0$ and $M$ contains all roots of 1, then for any $n>1$, the field $F$ has an extension of degree $n$ satisfying conditions (A) - (C).
\end{example}


In the notation above, we have the following


\begin{prop} \label{p:one extension}
Let $c\in F^* \backslash {F^*}^2$. Then
there exists a regular field extension $F_c/F$ such that

(1) the homomorphism $\res_{F_c/F} : \Br(F) \longrightarrow \Br(F_c)$ is injective;

(2) the field $F_c(\sqrt{c})$ splits the algebras $\cA_{F_c}$ and $\cB_{F_c}$.

Moreover, if a field extension $K/F$ satisfies conditions (A) - (C) and $c \not \in {K^*}^2$, then

(3) there is no $a\in F_c$ such that $da \in {F_c K^*}^2$;

(4) the composite $F_c K$ does not split $\cB_{F_c}$;

(5) the composite $F_c K(\sqrt{d})$ splits $\cB_{F_c}$ but does not split $\cA_{F_c}$.

\end{prop}


\noindent {\it {Proof}}.
Let $F(x)$ be a purely transcendental extension of $F$ of transcendence degree 1. Let also
$$
\cC := \cA_{F(x)}\otimes (c,x)
$$
be a biquaternion $F(x)$-algebra
and $F_1$ the function field of the Severi-Brauer variety of the algebra $\cC$.

Now let $F_1(y)$ be a purely transcendental extension of $F_1$ of transcendence degree 1 and
$$
\cD := \cB_{F_1(y)}\otimes (c,y)
$$
a biquaternion $F_1(y)$-algebra.
Let also $F_c$ be the function field of the Severi-Brauer variety of the algebra $\cD$.

Since the kernel of the restriction homomorphism $\res_{F_1/F(x)}$ is generated by the class of the algebra $\cC$, then the homomorphism $\res_{F_1/F}$ is injective. Indeed,
$\cC$ ramifies at the discrete valuation (trivial on $F$) of $F(x)$ defined by the polynomial $x$. Hence $[\cC] \ne [\cQ_{F(x)}]$ for any central simple $F$-algebra $\cQ$.

Note that $F_1$ splits $\cC$. Then $\cA_{F_1} \cong (c,x)_{F_1}$ and $F_1(\sqrt{c})$ splits $\cA_{F_1}$.

Let $K/F$ be a field extension satisfying conditions (A) - (C). Since $F_1/F$ is a regular field extension, then there
is no $a\in F_1$ such that $da \in {F_1 K^*}^2$. In particular, this means that $dc \not \in {F_1 K^*}^2$.
Moreover, $F_1 K/K$ is a regular extension of $K$. Thus, if $c \not \in {K^*}^2$, then $c \not \in {F_1 K^*}^2$.

The composite $F_1K (\sqrt{d})$ is the function field of the Severi-Brauer variety of the $K(\sqrt{d})(x)$-algebra $\cC_{K(\sqrt{d})(x)}$.  Hence the kernel of the
restriction homomorphism $\res_{F_1 K(\sqrt{d})/K(\sqrt{d})(x)}$ is generated by the class of $\cC_{K(\sqrt{d})(x)}$. Since $dc \not \in {K^*}^2$, then
$c \not \in K(\sqrt{d})^{*2}$ and $\cC_{K(\sqrt{d})(x)}$
ramifies at the discrete valuation (trivial on $K(\sqrt{d})$) of $K(\sqrt{d})(x)$ defined by the polynomial $x$, but $\cA_{K(\sqrt{d})(x)}$ is unramified at this valuation,
hence $[\cA_{K(\sqrt{d})(x)}]\ne [\cC_{K(\sqrt{d})(x)}]$ and $F_1 K(\sqrt{d})$ does not split $\cA_{F_1}$. Analogously,  the composite
$F_1 K$ does not split $\cB_{F_1}$.

Thus the extension $F_1 K/F_1$
satisfies conditions (A)-(C) with respect to the algebras $\cA_{F_1}$ and $\cB_{F_1}$.

The field $F_c$ satisfies conditions (1)-(5) of the proposition by the same arguments as for the field $F_1$. We just replace the ground field $F$ by $F_1$ and the extension
$K/F$ by $F_1 K/F_1$.
\qed


\begin{remark}
Conditions (1) and (3)-(5) of Proposition \ref{p:one extension} say that if $K/F$ is a field extension satisfying conditions (A)-(C), then the extension $F_c K/F_c$
satisfies conditions (A)-(C) with respect to the algebras $\cA_{F_c}$ and $\cB_{F_c}$.
\end{remark}

\medskip

In the notation above, we also have the following

\begin{proposition} \label{pr:one field}
Let $U :=\{c \in F^* \backslash {F^*}^2   \mid F(\sqrt{c}) \mbox{ splits }  \cA  \mbox{ or }  \cB  \}$.
There exists a regular field extension $E(F)/F$ such that

(1) the homomorphism $\res_{E(F)/F} : \Br(F) \longrightarrow \Br(E(F))$ is injective;

(2) the field $E(F)(\sqrt{c})$ splits the algebras $\cA_{E(F)}$ and $\cB_{E(F)}$ for any $c \in U$.

Moreover, if a field extension $K/F$ satisfies conditions (A) - (C), then

(3) there is no $a\in E(F)$ such that $da \in {E(F) K^*}^2$;

(4) the composite $E(F) K$ does not split $\cB_{E(F)}$;

(5) the composite $E(F) K(\sqrt{d})$ splits $\cB_{E(F)}$ but does not split $\cA_{E(F)}$.

\end{proposition}

\noindent {\it {Proof}}.
Note that for any field extension $K/F$ satisfying conditions (A) - (C), $U \cap {K^*}^2=\emptyset$ since $K$ does not split $\cA$ and $\cB$.

Let $<$ be a well-ordering on the set $U$ and let $c_0$ denote its least element.
Set $E^{c_0} := F_{c_0}$, where the field $F_{c_0}$ is constructed in Proposition \ref{p:one extension}.

For $c\in U$, $c \ne c_0$, set
$$
E^{<c} := \bigcup_{c' < c } E^{c'} \mbox{ and }
E^c : ={E^{<c}_c},
$$
where the field $E^c$ is obtained by applying Proposition \ref{p:one extension} to the field
$E^{<c}$ and the element $c \in E^{<c}$ and the algebras $\cA_{E^{<c}}$ and $\cB_{E^{<c}}$.
Define also
$E(F) := \bigcup_{c\in U} E^c$.

By Proposition \ref{p:one extension} and transfinite induction, the field $E(F)$ satisfies conditions (1)-(5) of the proposition.
\qed


\medskip


\begin{theorem}  \label{th: main}
Let $F$ be a field such that there are two non-isomorphic quaternion $F$-algebras $\cA$ and $\cB$.
There exists a regular field extension $L/F$ with the following properties:

(1) $\cA_L$ and $\cB_L$ are division algebras and ${\bf gen} (\cA_L) = {\bf gen} (\cB_L)$;

(2) If $K/F$ is a field extension of degree $n$ satisfying properties (A) - (C) with respect to the algebras $\cA$ and $\cB$, then the matrix algebras
$M_n(\cA_L)$ and $M_n(\cB_L)$ do not have the same maximal subfields;

(3) If $K$ is a field from the previous item and $C$ is a central division $F$-algebra of exponent $n$ which is split by $K$, then $C_L$ is a division algebra
of exponent $n$ and the algebras
$\cA_L \otimes \cC_L$ and $\cB_L \otimes \cC_L$ do not have the same maximal subfields.
\end{theorem}


\noindent {\it {Proof}}.
Let $K_0 := F$. We recursively define $K_i$, $i\in {\mathbb{Z}}_ { > 0}$, to be the field
$E(K_{i-1})$ constructed by applying
Proposition \ref{pr:one field} to the field $K_{i-1}$ and the algebras $\cA_{K_{i-1}}$ and $\cB_{K_{i-1}}$.
Let also $L := \bigcup_{i\ge 0} K_i$.

By induction and Proposition \ref{pr:one field},
the homomorphism $\res_{L/F} : \Br(F) \longrightarrow \Br(L)$ is injective. Hence
$\cA_L$ and $\cB_L$ are non-isomorphic division algebras.

Assume that $M$ is a maximal subfield of $\cA_L$.
Then there exists $i\ge 0$ such that
$M= L M'$, where $M'$ is a quadratic extension of $K_i$ that splits $\cA_{K_i}$.
By the construction of $K_{i+1}$, the field
$K_{i+1} M'$ splits the algebra $\cB_{K_{i+1}}$. Hence
$M$ splits $\cB_L$. Analogously, every maximal subfield of $\cB_L$ splits $\cA_L$.
Thus the algebras $\cA_L$ and $\cB_L$ have the same family of maximal subfields, i.e.,
${\bf gen}(\cA_L) = {\bf gen}(\cB_L)$.

Assume that $K/F$ is
a field extension of degree $n$ satisfying conditions (A) - (C) with respect to the algebras $\cA$ and $\cB$.
By induction and Proposition \ref{pr:one field},
the composite $L K(\sqrt{d})$ splits $\cB_L$ but does not split $\cA_L$.
Then $L K(\sqrt{d})$
embeds in $M_n(\cB_L)$ but does not embed in $M_n(\cA_L)$. Hence $M_n(\cA_L)$ and $M_n(\cB_L)$ do not have the same maximal subfields.

Finally, let $\cC$ be a central division $F$-algebra of exponent $n$ which is split by $K$.
Since the homomorphism $\res_{L/F}$ is injective,
then the exponent of $\cC_L$ is $n$. Since the exponent of $\cC_L$ divides its index, then $\cC_L$ is a division algebra.
The composite $L K(\sqrt{d})$
splits $\cB_L\otimes \cC_L$ but does not split $\cA_L\otimes \cC_L$. This means that $\cA_L \otimes \cC_L$ and $\cB_L \otimes \cC_L$ do not have the same maximal subfields.
\qed

\medskip


\begin{corollary} \label{cor:matrix algebras}
There exists a field $L$ such that there are two quaternion division $L$-algebras
$\cD$ and $\cD'$ such that ${\bf gen} (\cD) = {\bf gen} (\cD')$, but for any $n>1$
the matrix algebras $M_l(\cD)$ and $M_l(\cD')$ do not have
the same maximal subfields.
\end{corollary}


\noindent {\it {Proof}}. Let $F$ be a field such that there are two non-isomorphic quaternion $F$-algebras $\cA$ and $\cB$ and
for any $n>1$, the field $F$ has an extension of degree $n$ satisfying conditions (A) - (C) with respect to
the algebras $\cA$ and $\cB$. Let $L$ be the field constructed in Theorem \ref{th: main}. By Theorem \ref{th: main}, the algebras $\cD:=\cA_L$ and $\cD':=\cB_L$
have the required properties.

\qed


\end{document}